\documentclass[12pt]{article}
\usepackage[utf8]{inputenc}
\usepackage{amsfonts}
\usepackage{amssymb}
\usepackage{amsthm}
\usepackage[fleqn]{amsmath}
\usepackage{epsf, epsfig}
\usepackage{amssymb,latexsym}
\usepackage{color, lineno}

\usepackage{subfigure}
\usepackage{graphicx}
\usepackage{lipsum}
\usepackage{pinlabel}
\usepackage{epstopdf}

\newtheorem{lem}{Lemma}[section]
\newtheorem{theo}{Theorem}[section]

\newtheorem{cor}{Corollary}[section]
\newtheorem{con}{Conjecture}[section]

\newtheorem{prob}{Problem}[section]

\def \iff {if and only if }

\newcommand{\beeq}{\begin{equation}}
\newcommand{\eneq}{\end{equation}}

\newcommand{\beeqn}{\begin{eqnarray*}}
\newcommand{\eneqn}{\end{eqnarray*}}

\def \realr {{\cal GR}} 
\def \intr {{\cal GI}}

\newcommand {\relabel}[1] {\label{#1} \red{[*: #1]}}\newcommand {\rebibitem}[1] {\bibitem{#1} \red{[*: #1]}} 
\def\relabel {\label} \def\rebibitem {\bibitem}  

\begin{document}

\title{A survey on the study of real zeros of flow polynomials\thanks{This paper was partially supported 
by NTU AcRF project (RP 3/16 DFM) of Singapore).}}
\author{
Fengming Dong\thanks{Corresponding author.
Email: fengming.dong@nie.edu.sg}\\
\small Mathematics and Mathematics Education\\
\small National Institute of Education\\
\small Nanyang Technological University, Singapore 637616\\
}
\date{}

\maketitle

\begin{abstract}
For a bridgeless graph $G$, its flow polynomial 
is defined to be the function $F(G,q)$  which counts 
the number of nonwhere-zero $\Gamma$-flows
on an orientation of $G$
whenever $q$ is a positive integer
and $\Gamma$ is an additive Abelian group of order $q$.
It was introduced by Tutte in 1950
and the locations of zeros of this polynomial
have been studied by many researchers. 
This article gives a survey on the results 
and problems on the study of real zeros 
of flow polynomials. 
\end{abstract}

Keywords: bridgeless graph, nowhere-zero flow, flow polynomial,  flow root, zero-free interval

\section{Introduction}\relabel{intro}

\def \setz {\mathbb Z}
\def \seti {{\cal I}}

Let $G=(V,E)$ be a finite multigraph with vertex set 
$V$ and edge set $E$ and let $D$ be an orientation of $G$. 
For any additive Abelian group $\Gamma$, 
a {\it $\Gamma$-flow} on $D$ is a mapping 
$\phi: E\rightarrow \Gamma$ such that 
\begin{equation}\relabel{sec1-eq1}
\sum_{e\in A^+(v)}\phi(e)=\sum_{e\in A^-(v)}\phi(e)
\end{equation}
holds for every vertex $v$ in $G$,
where 
$A^+(v)$ (resp. $A^-(v)$) is the set of arcs in $D$ 
with tail $v$
(resp. the set of arcs in $D$ with head $v$).
A $\Gamma$-flow $\phi$ on $D$
is called a {\it nowhere-zero}
$\Gamma$-flow on $D$ if $\phi(e)\ne 0$ holds 
for all $e\in E$.

By applying (\ref{sec1-eq1}),
one can shows 
that $\sum\limits_{e\in A^+(S)}\phi(e)=
\sum\limits_{e\in A^-(S)}\phi(e)$ holds for 
any $S\subseteq V$ and any $\Gamma$-flow $\phi$ on $D$,
where  $A^+(S)$ (resp. $A^-(S)$) is 
the set of non-loop arcs in $D$ 
with tails in $S$ and heads in $V-S$
(resp. the set of non-loop arcs in $D$ 
with heads in $S$ and tails in $V-S$). 
Thus
there is no nowhere-zero $\Gamma$-flow on $D$
when $G$ contains a bridge.

For any integer $q\ge 2$, 
a {\it nowhere-zero $q$-flow} on $D$ is defined to be 
a nowhere-zero ${\mathbb Z}$-flow $\psi$ on $D$ 
such that $|\psi(e)|\le q-1$ holds for all $e\in E$,
where ${\mathbb Z}$ is the additive group 
consisting of integers. 
Due to a result of Tutte~\cite{tut2}, 
these two kinds of nowhere-zero flows have the same 
property of existence. 

\begin{theo}[\cite{tut2}]\relabel{tut-th0}
For any orientation $D$ of $G$ and any positive integer $q$,
there exists a nowhere-zero $q$-flow on $D$ 
if and only if there exists a nowhere-zero $\Gamma$-flow on $D$,
where $q$ is the order of $\Gamma$.
\end{theo}


The flow polynomial of a graph was introduced by 
Tutte~\cite{tut2} in 1950. 
For any positive integer $q$, let $F(G,q)$ be the 
the number of nowhere-zero $\Gamma$-flows on $D$,
where $\Gamma$ is an additive Abelian group of order $q$.
It is not difficult to verify that 
the definition of $F(G,q)$ does not depend on 
the selections of $D$ and $\Gamma$ but on $G$ and $q$.
For any positive integer $q$, 
Theorem~\ref{tut-th0} implies that 
$F(G,q)>0$ \iff there exists a nowhere-zero $q$-flow
on $D$ for any orientation $D$ of $G$.

Note that $F(G,q)$ is not equal to 
the number of nowhere-zero $q$-flows. 
For example, if $G$ is the graph with one vertex and one loop,
then $F(G,q)=q-1$ while the number of nowhere-zero $q$-flows
is $2(q-1)$. 

The function $F(G,q)$ can also be determined by 
the following properties \cite{tut}: 
\begin{equation}\relabel{sec1-eq2}
F(G,q)=
\left \{
\begin{array}{ll}
1, &\mbox{if }E=\emptyset;\\
0, &\mbox{if }G \mbox{ has a bridge};\\
F(G_1,q)\cdots F(G_k,q), &\mbox{if }G=G_1\cup \cdots \cup G_k;\\
(q-1)F(G-e,q), &\mbox{if }e \mbox{ is a loop in }G;\\
F(G/e, q) -F(G- e, q),  &\mbox{if }
e \mbox{ is neither a loop nor a bridge,}
\end{array}
\right.
\end{equation}
where $G/e$ and $G - e$ are the graphs obtained from $G$ by contracting $e$ and deleting $e$ respectively
and $G_1\cup \cdots \cup G_{k-1} \cup G_k$
is the disjoint union of graphs 
$G_1, G_2, \cdots, G_{k-1}$ and $G_k$. 
By applying the properties in (\ref{sec1-eq2}), 
$F(G,q)$ can be expressed in terms of 
the Tutte polynomial $T_G(x,y)$ of $G$:
\begin{equation}\relabel{sec1-eq3}
F(G,q)=\sum_{E'\subseteq E}(-1)^{|E|-|E'|}
q^{|E'|+c(E')-|V|}
=(-1)^{|E|-|V|+c(E)}T_G(0, 1-q),
\end{equation}
where 
$c(E')$ is the number of components of 
the spanning subgraph $(V,E')$ of $G$ 
and
\begin{equation}\relabel{sec1-eq4}
T_G(x,y)=\sum_{E'\subseteq E}
(x-1)^{c(E')-c(E)}(y-1)^{|E'|-|V|+c(E')}.
\end{equation}
Both (\ref{sec1-eq2}) and (\ref{sec1-eq3}) show that 
$F(G,q)$ is a polynomial in $q$.
Thus the variable $q$ in $F(G,q)$ 
can be considered as a real or complex number 
for the study of its algebraic or other properties. 

The zeros of the flow polynomial $F(G,q)$ 
are called {\it the flow roots} of $G$. 
This article focuses on giving a review on 
the study of real flow roots of graphs. 

In Section~\ref{sect-basic}, we will introduce some basic 
results on flow polynomials.
Because the flow roots of any plane graph 
are exactly the non-zero roots of 
the chromatic polynomial of its dual plane graph
(see Theorem~\ref{Tut-th1}),
we will give a short review on 
the study of real roots of chromatic polynomials
in Section~\ref{sect-chro}.
The other sections are arranged below: 
\begin{enumerate}
\item 
(Section~\ref{zero-free}) 
determining maximal zero-free intervals of flow polynomials
of the form $(1,a)$ for $1<a\le 2$
and in particular, 
searching for graphs having no flow roots in 
the interval $(1,2)$;
\item 
(Section~\ref{zero-free2}) 
searching for near-cubic graphs
which have no flow roots in the interval $(2,3)$;
\item 
(Section~\ref{roots2})
the multiplicity of flow root ``$2$" 
for near-cubic graphs;
\item 
(Section~\ref{roots-ge4})
the existence of flow roots larger than $4$;
\item (Section~\ref{real zero})
the existence of graphs 
which have real flow roots only 
but also contain non-integral real flow roots. 
\end{enumerate}

\section{Basic properties of flow polynomials
\relabel{sect-basic}}

A graph $G=(V,E)$ is said to be 
{\it non-separable} if 
either $|E|\le |V|=1$ or 
$G$ is connected without loops or cut-vertices, 
where a vertex $x$ in $G$ is called a 
{\it cut-vertex} if $G-x$ has more components than $G$ has.
We say $G$ is {\it separable} if it is not non-separable. 

A {\it block} of $G$ is 
a maximal subgraph of $G$ 
with the property that it is non-separable.
Clearly, if $|E|+|V|\ge 3$, then each loop of $G$ 
is considered as a block. 
A block is said to be {\it trivial} if 
its order is $1$ and its size is $0$. 
So a trivial block is an isolated vertex in $G$.
By (\ref{sec1-eq2}), $F(G,q)=F(G-V_0,q)$ holds,
where $V_0$ is the set of isolated vertices in $G$.
Thus we may assume that $G$ has no isolated vertices 
when we study $F(G,q)$.

By (\ref{sec1-eq2}), $F(G,q)=0$ whenever $G$ contains bridges.
Also by (\ref{sec1-eq2}), if $G$ is separable, 
then $F(G,q)$ can be expressed as the multiplication 
of flow polynomials of its components or blocks.

\begin{lem}
\relabel{block-factor}
If $G_1, G_2, \cdots, G_k$ are the components of $G$
or the blocks of a connected graph $G$, then 
\beeq
F(G,q)=\prod_{1\le i\le k}F(G_i,q).
\eneq
\end{lem}

For any $S\subseteq E$, $S$ is called 
an {\it edge-cut} of $G$ if $S$ is the set 
of edges with ends  in $V_1$ 
and in $V_2$ respectively, where 
$\{V_1,V_2\}$ is a partition of $V$
with $V_i\ne \emptyset$ for $i=1,2$.
Such an edge-cut $S$ is also written as $S=(V_1,V_2)$.
For any non-separable graph $G$, 
$F(G,q)$ can be further factorized in any one of the 
following cases, due to Jackson~\cite{jac3} (see \cite{dong2, jac2,jac4} also):
\begin{enumerate}
\item $G-e$ is separable for some edge $e$;
\item $G$ has a proper edge-cut $S$ with $2\le |S|\le 3$,
where an edge-cut $S$ is said to be {\it proper} if $G-S$ has no isolated vertices.
\end{enumerate}
For any graph $G$ and any two vertices $u$ and $v$ in $G$,
let $G+uv$ denote the graph obtained by adding a new edge 
joining $u$ and $v$.

\begin{lem}[\cite{jac3}]
\relabel{v-edge}
Let $G$ be a bridgeless connected graph
and $e$ be an edge of $G$ joining $u_1$ and $u_2$.
If $G-e$ is separable and 
$H_1$ and $H_2$ are edge-disjoint subgraphs of 
$G-e$ with $E(H_1) \cup E(H_2) = E(G - e)$, 
$V(H_1)\cup V(H_2) =V(G)$,
$V(H_1)\cap V(H_2) = \{v\}$ for some vertex $v$ of $G$
and $u_i\in V (H_i)$ for $i=1,2$, as shown in 
Figure~\ref{p55-f1}, then 
\beeq
F(G, q) =
\frac{F(G_1, q)F(G_2, q)}{q -1},
\eneq
where $G_i=H_i+vu_i$ for $i=1,2$.
\end{lem}

\begin{figure}[h!]
\centering 
\includegraphics[width=0.6\textwidth]{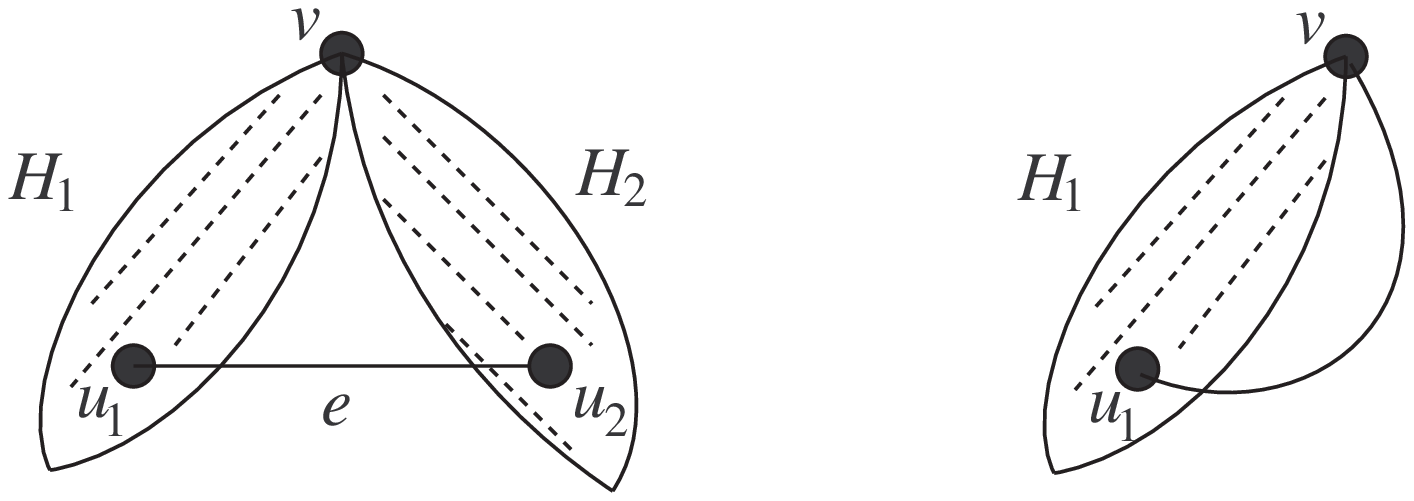}


$G$ \hspace {4.5 cm} $G_1$

\caption{\relabel{p55-f1} $G-e$ is separable.}
\end{figure}


\begin{lem}[\cite{jac3}]
\relabel{2-edge}
Let $S=(V_1,V_2)$ be an edge-cut of 
a $2$-edge connected graph $G$  
and let $H_i=G[V_i]$ be the subgraph of $G$ 
induced by $V_i$ for $i=1,2$,
as shown in 
Figure~\ref{p55-f2} when $|S|=2$. 
For $i=1,2$, 
let $G_i$ be obtained from $G$ by contracting $E(H_{3-i})$. 
If $2\le |S|\le 3$, then 
\beeq\relabel{2-block-eq1}
F(G, q) =
\frac{F(G_1, q)F(G_2, q)}
{(q -1)_{|S|-1}},
\eneq
where $(x)_k$ is the polynomial $x(x-1)\cdots (x-k+1)$. 
\end{lem}

\begin{figure}[h!]
\centering 
\includegraphics[width=0.55\textwidth]{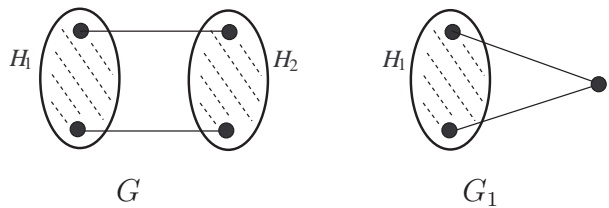}

$G$ \hspace {4 cm} $G_1$

\caption{\relabel{p55-f2} $G$ has a $2$-edge-cut.}
\end{figure}

\noindent {\bf Remark}:  
If $S$ is a non-proper 
edge-cut with $2\le |S|\le 3$,
then $V_i$ contains one vertex only for some $i$, 
implying that $G_{3-i}\cong G$,
where $G_i$ is a graph defined in Lemma~\ref{2-edge}.
Thus, only when $S$ is a proper edge-cut of $G$
with $2\le |S|\le 3$, 
Lemma~\ref{2-edge} can be applied to 
express $F(G,q)$ in terms of multiplications 
of flow polynomials of graphs with smaller orders.  

For any edge $e$ of $G$, if one end of $e$ is of degree $2$,
then the results of (\ref{sec1-eq2}) imply 
that $F(G,q)=F(G/e,q)$.
Thus, by Lemmas~\ref{block-factor},~\ref{v-edge}
and~\ref{2-edge},
the study of zeros of flow polynomials can be focused on 
non-separable and $3$-edge connected graphs 
which do not contain any proper 3-edge-cut
and do not contain any edge whose removal results in 
a separable graph.


\section{Real zeros of chromatic polynomials
\relabel{sect-chro}}

The flow polynomial $F(G,q)$ is considered as the dual polynomial 
of {the chromatic polynomial} $P(G,q)$,
mainly due to their close relation given in 
Theorem~\ref{Tut-th1},
where 
$P(G,q)$ is defined to be the function which counts 
the number of proper $q$-colourings of $G$
whenever $q$ is a positive integer. 
Clearly, for a positive integer $q$, 
$G$ admits a proper $q$-colouring 
\iff $P(G,q)>0$.

\begin{theo}[\cite{tut3}] \relabel{Tut-th1}
$P(G,q)=qF(G^*,q)$ holds for any connected plane graph $G$, 
where $G^*$ is the dual plane graph of $G$.
\end{theo}

The chromatic polynomial $P(G,q)$ was introduced by 
Birkhoff~\cite{bir0} in 1912 in the hope of proving 
the {\it four-color theorem}
(i.e., $P(G,4)>0$ holds for any non-loop planar graph $G$).
This function $P(G,q)$ is indeed a polynomial in $q$, 
as its definition implies that for any integer $q\ge 1$, 
\begin{equation}\relabel{sec1-eq0}
P(G,q)=\sum_{1\le k\le |V(G)|} \alpha_k(G) (q)_k,
\end{equation}
where 
$\alpha_k(G)$ is the number of partitions of 
$V(G)$ into exactly $k$ non-empty independent sets. 
Thus the variable $q$ in the function $P(G,q)$ 
can be considered as a complex number. 
If $P(G,q)=0$, then $q$ is called a {\it chromatic root} of $G$. 

By Theorem~\ref{Tut-th1},
the flow roots and chromatic roots of planar graphs 
have the same distribution, except that 
`$0$' is a chromatic root of every graph 
but not a flow root of any graph. 

Due to Sokal's result below on complex zeros of 
chromatic polynomials and Tutte's result in 
Theorem~\ref{Tut-th1}, 
the complex flow roots of 
planar graphs are dense everywhere in the
whole complex plane with the possible 
exception of the disc $|z - 1| < 1$.

\begin{theo}[\cite{sok2004}]\relabel{sok-th1}
The complex zeros of chromatic polynomials of 
planar graphs are dense everywhere
in the whole complex plane 
with the possible exception of the disc $|z-1| < 1$.
\end{theo}

This section focuses on giving a review on 
the study of real chromatic roots of graphs, 
and in particular, of planar graphs. 
For other results or problems on chromatic polynomials, 
the reader can refer to \cite{bir, dong3, rea, read1988}.

For general graphs, 
$(-\infty,0), (0,1)$ and $(1,32/27]$ are 
the only maximal zero-free intervals for 
all chromatic polynomials,
where an interval is said to be {\it zero-free} 
for a function if it has no zero in this interval.

The first two zero-free intervals for all chromatic 
polynomials follow 
directly from parts (i) and (iii) of the following result.

\begin{theo}[\cite{rea, read1988}]\relabel{zero-free1}
Let $G$ be a non-loop graph of order $n$
and component number $c$. Then 
\begin{enumerate}
\item $(-1)^n P(G,q)>0$ for all real $q<0$; 
\item the multiplicity of the root $0$
of $P(G,q)$ 
is equal to $c$; 
\item $(-1)^{n-c} P(G,q)>0$ for all real $0<q<1$;
\item the multiplicity of the root $1$ 
of $P(G,q)$ 
is equal to the number of non-trivial blocks.
\end{enumerate}
\end{theo}

The third such zero-free interval $(1,32/27]$
was  due to Jackson~\cite{jac1}.

\begin{theo}[\cite{jac1}]\relabel{zero-free3}
Let $G$ be a non-loop connected graph of order $n$ 
and block number $b$. Then 
$(-1)^{n+b-1}P(G, q) >0$ 
holds for all real $q$ in $(1,32/27]$.
\end{theo}

Thomassen~\cite{tho} showed that 
these three intervals are the only 
zero-free intervals for all chromatic polynomials.

\begin{theo}[\cite{tho}]\relabel{zero-free4}
For any interval $(a,b)$ with $32/27\le a<b$, 
there exists a graph which has chromatic roots 
in the interval $(a,b)$.
\end{theo}


Now we focus on planar graphs. 
By a combinatorial approach, it is not difficult to 
prove that every non-loop 
planar graph admits a proper 5-colouring
(i.e., $P(G,5)>0$ holds for any non-loop planar graph $G$).
This result can also be proved by induction 
with an algebraic approach.
Actually, Birkhoff and Lewis \cite{bir} 
proved 
that 
$P(G,q)>0$ holds for all non-loop planar graphs $G$ 
and all real numbers $q\ge 5$.
Thus, $[5,\infty)$ is a zero-free interval 
for chromatic polynomials of all non-loop planar graphs. 


\begin{theo}[\cite{bir}]\relabel{bir-th1}
$P(G,q)>0$ holds for all non-loop planar graphs $G$
and all real numbers $q$ in the interval 
$[5,\infty)$.
\end{theo}

Note that Theorem~\ref{bir-th1} does not hold if $G$ 
is not restricted to planar graphs as 
$P(G,q)=0$ whenever $q$ is an integer with
$0\le q <\chi(G)$.

Birkhoff and Lewis \cite{bir} also conjectured that 
$[4,5)$ is a zero-free interval 
for chromatic polynomials of all non-loop planar graphs.

\begin{con}[\cite{bir}]\relabel{bir-con1}
$P(G,q)>0$ holds for all non-loop plane graphs $G$
and all real numbers $q$ in the interval 
$[4,5)$.
\end{con}

Conjecture~\ref{bir-con1} includes 
the four-color conjecture (i.e., the case $q=4$).
The four-color conjecture was first proved  
Appel and Haken~\cite{appel1977} in 1977 and was reproved by 
Robertson, Sanders, Seymour and Thomas~\cite{rob1997} in 1997.
However, the study of Conjecture~\ref{bir-con1} 
has no other progress
and it is even unknown if there exists 
a real number $\epsilon>0$ such that 
$P(G,q)>0$ holds 
for all plane graphs $G$ and 
all real $q\in (5-\epsilon,5)$.  

The number $\tau=\sqrt 5/2+1/2\approx 1.618033$ is called the 
{\it golden ratio}
which is the positive real root of the equation $q^2=q+1$. 
The following results   
due to Tutte~\cite{tut02, tut01},
Thomassen~\cite{tho} 
and Perret and Thomassen~\cite{per 2016} 
present some interesting relations between $\tau$ 
and chromatic roots of planar graphs.

\begin{theo}
\begin{enumerate}
\item \cite{tut02, tut01}
For any plane triangulation $G$ with $n$ vertices, 
$$
0<|P(G,\tau+1)|\le \tau^{5-n} 
$$
and 
$$
P(G,\tau+2)=(\tau+2) \tau^{3n-10}(P(G,\tau+1))^2>0;
$$
\item \cite{per 2016, tho} chromatic roots of planar graphs
are dense everywhere in the interval $[32/27,4)$
with the possible exception\label{r2-q5}
of a small interval $(t_1,t_2)$ around 
the number $\tau+2\approx 3.618033$,
where $t_1\approx 3.618032$ and $t_2 \approx 3.618356$.
\end{enumerate}\relabel{zero-free5}
\end{theo}

Note that in Theorem~\ref{zero-free5} (i), 
$P(G,\tau+1)\ne 0$. 
Actually, $\tau+1=(3+\sqrt 5)/2$ 
cannot be a chromatic root of any 
graph $G$. Otherwise, $(3-\sqrt 5)/2\approx 0.381966$ 
is also a chromatic root of $G$, contradicting 
Theorem~\ref{zero-free1} (iii).
More generally, for any positive rational numbers 
$a$ and $b$, if $\sqrt b$ is not rational and 
$a-1<\sqrt b$, then $a+\sqrt b$ is not a chromatic root
of any graph.

\section{Zero-free intervals within $(1,2)$
for flow polynomials}
\relabel{zero-free}

Let  $G =(V, E)$ be a bridgeless connected graph. 
Applying (\ref{sec1-eq2}), 
it is not difficult to show that 
$(-1)^{|E|-|V|+1}F(G, q)>0$ holds for all real $q$ 
in the interval $(-\infty,1)$.
But $F(G,1)=0$ if $E\ne \emptyset$.
Thus $(-\infty,1)$ is a maximal zero-free 
interval of all flow polynomials.
This is part of the results below due to Wakelin~\cite{wak}.

\begin{theo}[\cite{wak}]\relabel{wak-th}
Let $G =(V, E)$ be a bridgeless connected graph
with block number $b$. Then
\begin{enumerate}
\item $(-1)^{|E|-|V|+1}F(G, q)>0$ holds for all real $q$ 
in $(-\infty,1)$;
\item $F(G, q)$ has a zero of multiplicity $b$ at $q =1$;
\item $(-1)^{|E|-|V|+b+1}F(G, q)>0$ holds for all real $q$ 
in $(1, 32/27]$.
\end{enumerate}
\end{theo}

By Theorem~\ref{wak-th}, $(1, 32/27]$ is the second
zero-free interval for all flow polynomials. 
By Theorems~\ref{Tut-th1} and~\ref{zero-free5} (ii), 
this zero-free interval $(1, 32/27]$ for flow 
polynomials is also maximal,
and 
there is no other 
zero-free interval $(a,b)$ 
for all flow polynomials with $32/27<a<b\le 4$, 
unless $t_1\le a<b\le t_2$,
where $t_1,t_2$ are number stated 
in Theorem~\ref{zero-free5} (ii).
But it is unknown if there exists a 
zero-free interval $(a,b)$ 
for all flow polynomials with $4<a<b$. 
The study of this problem will be reviewed in 
Section~\ref{roots-ge4}.

Now we consider maximal zero-free intervals
of flow polynomials of some subsets of graphs. 
A {\it plane near-triangulation} is a non-loop connected plane graph in which at most one face is not bounded by a cycle of order 3. 
Birkhoff and Lewis~\cite{bir} proved that 
any near-triangulation does not have 
any chromatic root in $(1,2)$. 

\begin{theo}[\cite{bir}]\relabel{bir-th2}
If $G$ is a plane near-triangulation, then 
$G$ has no chromatic roots in $(1,2)$. 
\end{theo}

By Theorems~\ref{Tut-th1} and~\ref{bir-th2}, 
we have the following equivalent statement.

\begin{theo}\relabel{bir-th3}
$(1,2)$ is a zero-free interval for flow polynomials of 
bridgeless planar graphs which have 
at most one vertex of degree not equal to $3$. 
\end{theo}

Jackson~\cite{jac3} showed that 
Theorems~\ref{bir-th3} holds no matter 
whether $G$ is planar or non-planar, 
as long as $G$ contains at most one vertex of degree 
larger than $3$.

\begin{theo}[\cite{jac3}]\relabel{jac-th1}
If $G$ is a bridgeless graph with at most one vertex of degree 
larger than $3$, then 
$G$ has no flow roots in $(1,2)$. 
\end{theo}

Theorem~\ref{jac-th1} was further generalized by 
Dong~\cite{dong1} who showed that the conclusion holds for 
every bridgeless  graph $G$ with $|W(G)| \le 2$,
where $W(G)$ is the set of vertices in $G$ 
which are of degrees larger than $3$.
But the conclusion fails for some bridgeless graphs $G$ with $|W(G)|=3$.
Actually, for any integer $k\ge 3$, it fails for 
some  bridgeless graphs $G$ with $|W(G)|=k$.
The graph shown in Figure~\ref{f1} 
has a flow root around $1.4301\cdots$, which is the only real zero
of the polynomial $q^3-5q^2+10q-7$.
This graph has the least size among all graphs 
which have flow roots in $(1,2)$.

\begin{figure}[h!]
\centering

\includegraphics[width=0.2\textwidth]{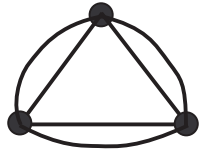}

\caption{\relabel{f1}
The graph with the smallest size and 
with a flow root in the interval $(1,2)$}
\end{figure}

For any integer $k\ge 0$, let $\Psi_k$ be the set of bridgeless connected graphs with $|W(G)|\le k$
and let $\xi_k$ be the supremum in $(1, 2]$ 
such that $(1, \xi_k)$ is a zero-free 
interval for flow polynomials of graphs in $\Psi_k$.
Clearly, $\xi_0, \xi_1, \xi_2,\cdots$ is a non-increasing 
sequence. 

For any $k\ge 2$, let 
$\Theta_k$ be the set of those graphs 
in $\Theta$ with exactly $k$ vertices, 
where $\Theta$ is the set of graphs defined by the two steps below:
\begin{enumerate}
\item[(i)] $Z_3\in \Theta$, where 
$Z_j$ is the graph with two vertices 
and $j$ parallel edges joining these two vertices; and 
\item[(ii)] $G(e)\in \Theta$ for every $G\in \Theta$ and 
every $e\in E(G)$, 
where $G(e)$ is the graph obtained from $G-e$ 
by adding a new vertex $w$ and for each end $u_i$ of $e$,
$i=1,2$, 
adding two parallel edges 
joining $w$ to $u_i$. 
\end{enumerate}

It is not difficult to check that 
$\Theta_k$ has exactly one graph for $k=2,3,4$,
and the graph shown in Figure~\ref{f1} is the only graph 
in $\Theta_3$. 

The value of $\xi_k$ can be determined by graphs 
in the finite set $\Theta_k$.

\begin{theo}[\cite{dong1}]\relabel{dong-th1}
For any integer $k\ge 2$, 
$\xi_k$ is the minimum value among the  
real flow roots in $(1,2)$ of all graphs in $\Theta_k$. 
\end{theo}

Applying Theorem~\ref{dong-th1}, 
the values of $\xi_k$'s for $k\le 5$ 
are determined. 

\begin{cor}[\cite{dong1}]\relabel{dong-cor1}
$\xi_k=2$ for $k=0,1,2$, 
$\xi_3=1.430\cdots$, 
$\xi_4=1.361\cdots$ and 
$\xi_5=1.317\cdots$, where 
the last three numbers 
are the real zeros of $q^3-5q^2+10q-7$,
$q^3-4q^2+8q-6$ 
and $q^3-6q^2+13q-9$
in $(1,2)$ respectively.
\end{cor}

Corollary~\ref{dong-cor1} tells that $(1,2)$ is a zero-free interval 
for flow polynomials of all bridgeless graphs $G$ with $|W(G)|\le 2$.
As $\xi_0, \xi_1, \xi_2, \xi_3, \cdots$ is non-increasing, 
Theorem~\ref{dong-th1} implies that for any $k\ge 3$, 
there exist graphs $G$ with $|W(G)|=k$ which have flow roots 
in $(1,2)$.
However, by the following result,
it is not true that 
every graph $G$ with $|W(G)|\ge 3$ 
contains flow roots in $(1,2)$.

For any $V_0\subseteq V(G)$, let $N_G(V_0)$ 
denote the set $\bigcup_{v\in V_0} N_G(v)$.

\begin{theo}[\cite{dong2}]\relabel{dong-th2}
For any bridgeless graph $G$,
if $G-W(G)$ has a component $G_0$ such that 
$W(G)\subseteq N_G(V(G_0))$ holds, then 
$(1,2)$ is a zero-free interval of $F(G,q)$.
\end{theo}

Theorem~\ref{dong-th1} and 
Corollary~\ref{dong-cor1} 
will be applied in the study of Problem~\ref{prob1}.

\section{Zero-free intervals within $(2,3)$}
\relabel{zero-free2}

For any vertex-disjoint graphs $G_1$ and $G_2$, 
let $G_1+G_2$ denote the graph with vertex set 
$V(G_1)\cup V(G_2)$ and edge set 
$E(G_1)\cup E(G_2)\cup \{xy: x\in V(G_1),y\in V(G_2)\}$.

Woodall~\cite{woo1992, woo} showed that every plane 
triangulation has no chromatic roots 
in the interval $(2, \alpha_1)$,
where $\alpha_1=2.546602\cdots$ is the unique real
chromatic root 
of the graph $C_4+\bar K_2$ in the interval $(2,3)$:
$$
P(C_4+\bar K_2,q)=q(q-1)(q-2)(q^3-9q^2+29q-32).
$$
Note that the dual graph of $C_4+\bar K_2$ is the
cube shown in Figure~\ref{f2}(a). 
By Theorem~\ref{Tut-th1}, $\alpha_1$ is 
the flow root of the cube
and  Woodall's result can be translated to 
the following one on flow polynomials.

\begin{figure}[h!]
\centering 

\includegraphics[width=0.6\textwidth]{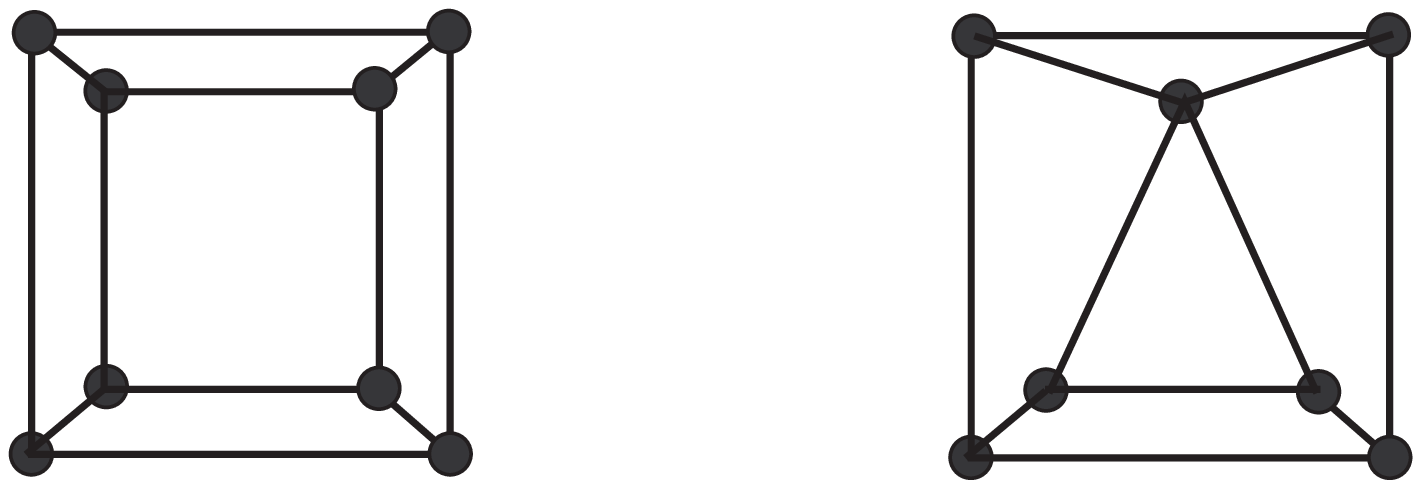}

\hspace{1.5 cm} (a) Cube \hspace{3 cm} (b) Contracted cube

\caption{\relabel{f2}
Cube and Contracted cube}
\end{figure}

\begin{theo}[\cite{woo1992, woo}]\relabel{woo-th1}
$(2,\alpha_1)$ is a zero-free interval 
for flow polynomials of bridgeless cubic planar graphs.
\end{theo}

Woodall~\cite{woo} conjectured that 
the conclusion of Theorem~\ref{woo-th1}
holds for the interval $(2.677814\cdots, 3)$,
where $2.677814\cdots$ is the unique real 
chromatic root
 of the graph $C_5+\bar K_2$ in the interval $(2,3)$:
$$
P(C_5+\bar K_2,q)=q(q-1)(q-2)(q-3)
(q^3-9q^2+30q-35).
$$

\begin{con}[\cite{woo}]\relabel{woo-con1}
$(2.677814\cdots, 3)$ is a zero-free interval for the flow polynomials of bridgeless cubic planar graphs. 
\end{con}

Jackson~\cite{jac4} extended Theorem~\ref{woo-th1}
to all cubic graphs. 

\begin{theo}[\cite{jac4}]\relabel{jac-th2}
$(2,\alpha_1)$ is a zero-free interval 
for flow polynomials of all bridgeless cubic graphs.
\end{theo}

We call $G$ a {\it near-cubic} graph if $|W(G)|\le 1$. 
Jackson~\cite{jac3} extended Theorem~\ref{jac-th1}
to near-cubic graphs for a smaller interval.  

\begin{theo}[\cite{jac3}]\relabel{jac-th3}
$(2,\alpha_2]$ \label{r1-q5}
is a zero-free interval 
for flow polynomials of bridgeless near-cubic graphs,
where $\alpha_2=2.225\cdots$ 
is the real zero in $(2, 3)$ of the polynomial 
$q^4 - 8q^3 + 22q^2 - 28q + 17$.
\end{theo}

Note that Theorem~\ref{jac-th3} does not hold 
if the number $\alpha_2$ is replaced by any larger number,
as Jackson \cite{jac3} has showed that 
for any $\epsilon>0$, there exists a near-cubic graph which 
has a flow root in $(\alpha_2,\alpha_2+\epsilon)$. 


Let $x\in V (G)$. 
A branch at $x$ is a path $P = xv_1v_2\cdots v_m$ such that $d_G(v_i) = 2$ for $1 \le i \le m - 1$ and 
$d_G(v_m)\ge 3$. 
We say that $v_1, v_2,\cdots, v_{m-1}$ 
are inner vertices of this branch at $x$. 
Let $G * x$ be the graph obtained from $G$ by deleting
$x$ and all the inner vertices of every branch at $x$. 

An edge-cut $S=(V_1,V_2)$ of a graph $G$ is said to be {\it cyclic} 
if $G[V_i]$ contains cycles for both $i=1,2$. 
A graph $G$ is said to be 
{\it cyclically $k$-edge-connected} 
if each cyclic edge-cut of $G$ has at least $k$ edges.
A graph is said to be {\it essentially $3$-connected} if it is a subdivision of a $3$-connected graph or the graph $Z_k$ for some $k\ge 3$.
For essentially $3$-connected near-cubic graphs, 
Jackson~\cite{jac4} obtained a zero-free interval
for their flow polynomials 
which is larger than the one in Theorem~\ref{jac-th3}.

\begin{theo}[\cite{jac4}]\relabel{jac-th4}
Let $G$ be a near-cubic graph 
with order $n$ and size $m$  
such that  
$G$ is essentially $3$-connected.
\begin{enumerate}
\item Suppose that for any $x\in V(G)$ with $d(x)\ge 5$,
\label{r2-q10}
$G * x$ is cyclically $3$-edge-connected.
\label{r1-q6}
Then $F(G, q)$ is nonzero with sign $(-1)^{m-n+1}$ 
for $q\in (2,\alpha_3)$, where $\alpha_3=2.43\cdots$ 
is a flow root of the contracted cube, 
shown in Figure~\ref{f2} (b),  whose flow polynomial 
is $(q-1)(q-2)(q^3-8q^2+23q-23)$.
\item 
Suppose that for any $x\in V(G)$ with $d(x)\ge 4$,
$G *x$ is cyclically 
$3$-edge-connected and $B(G*x)$ has at most
one component of order $1$,
where $B(G*x)$ is the subgraph of $G*x$ induced by 
the set of vertices of degrees $2$ in $G*x$.
Then $F(G, q)$ is non-zero with sign 
$(-1)^{m-n+1}$ for $q\in (2,\alpha_1)$.
\end{enumerate}
\end{theo}

\label{r2-q11}
It is mentioned in \cite{jac4} that 
$\tau^2=2.618\cdots$ is an accumulation
point of flow roots of cyclically 4-edge-connected cubic graphs,
as the Cartesian products $C_{2r} \times K_2$ 
for $r\ge 3$
are cyclically $4$-edge-connected and 
have flow roots tending to $\tau^2$ from
below as $r$ tends to infinity, where 
$\tau=\frac{1+\sqrt 5}{2}$ 
is the golden ratio. 
Jackson~\cite{jac4} proposed the following conjectures
on cyclically $4$-edge-connected cubic graphs. 

\begin{con}[\cite{jac4}]\relabel{con-jac2}
For all $\epsilon> 0$, there exist only 
finitely many cyclically $4$-edge-connected cubic
graphs with a flow root in $(2, \tau^2 -\epsilon)$.
\end{con}

\begin{con}[\cite{jac4}]\relabel{con-jac3}
Let $G$ be a cyclically $4$-edge-connected cubic graph. 
Then G has at most one flow root in $(2, \tau^2)$.
\end{con}

\section{Multiplicity of the flow root at 2 \relabel{roots2}}

For any bridgeless graph $G$ with $\delta(G)\ge 2$,
Theorem~\ref{wak-th} and 
the third equality in (\ref{sec1-eq2}) 
imply that 
the multiplicity of the flow 
root at $1$ of $G$ 
is equal to 
the total number of blocks in $G$.
This section focuses on the multiplicity of the flow root at $2$.

It is well known that a graph $G$ has a nowhere-zero 2-flow 
if and only if $G$ is an even graph (i.e., every vertex in $G$ 
is of an even degree). 
Thus Theorem~\ref{tut-th0} implies that 
$2$ is a flow root of $G$ \iff $G$ is not an even graph. 

Woodall~\cite{woo1992} showed that 
if $G$ is a 3-connected plane triangulation, 
then the multiplicity of chromatic roots of $G$ at $2$
is exactly $1$. 
By Theorem~\ref{Tut-th1}, Woodall's result is 
equivalent to that for any 
$3$-connected cubic planar graph $G$,
the multiplicity of flow roots of $G$ at $2$
is exactly $1$.
This result was extended to near-cubic graphs 
by Jackson~\cite{jac3}.

For any bridgeless graph $G$, let 
$$
q_2(G, q) =\frac{F(G, q)}{(q - 1)(q - 2)}.
$$

Jackson~\cite{jac3} showed that a non-separable near-cubic graph
$G$ with $|V(G)|\ge 2$
\label{r2-q12}
 is essentially 3-connected \iff 
the multiplicity of the flow root of $G$ at $2$ is at most $1$.

\begin{theo}[\cite{jac3}]\relabel{jac3-th4}
For any non-separable near-cubic graph $G$ with 
$|V(G)|\ge 2$, 
$q_2(G, q)$ is a polynomial in $q$. 
Furthermore:
\begin{enumerate}
\item if $G$ is not essentially $3$-connected, then 
$q_2(G, 2) = 0$;
\item if $G$ is essentially $3$-connected, then 
$q_2(G, 2)$ is non-zero with sign $(-1)^{m-n+1}$,
where $n$ and $m$ are the order and size of $G$.
\end{enumerate}
\end{theo}

\section{Flow roots larger than 4\relabel{roots-ge4}}

The Petersen graph has a flow root at $q=4$, as its flow polynomial is 
\begin{equation}\relabel{flow-petersen}
(q-1)(q-2)(q-3)(q-4)(q^2-5q+10).
\end{equation}
Thus, (\ref{flow-petersen})
and Theorem~\ref{tut-th0} 
imply that the Petersen graph does not admit 
a nowhere-zero $4$-flow.
By (\ref{sec1-eq2}), $F(G,q)=F(G/e,q)$ holds 
for any edge $e$ in $G$ which has one end of degree $2$.
Thus any subdivision of the Petersen graph 
does not admit a nowhere-zero $4$-flow.
Tutte~\cite{tut66} guessed that 
any bridgeless graph without a
subdivision of the Petersen graph
has a nowhere-zero $4$-flow.

\begin{con}[\cite{tut0}]\relabel{4-flow con}
Any bridgeless graph 
without a subdivision of the Petersen graph has 
a nowhere-zero $4$-flow.
\end{con}

By Theorem~\ref{Tut-th1},
this conjecture, called {\it Tutte's four-flow conjecture}, 
is obviously stronger 
than the Four Color Theorem,
as each plane graph does not have 
a subdivision of the Petersen graph.
For more details on Tutte's four-flow conjecture, 
the reader can refer to the surveys by Jaeger \cite{jae79}
and Younger \cite{you}. 

Tutte~\cite{tut0} also conjectured that 
any bridgeless graph admits a nowhere-zero $5$-flow,
known as {\it Tutte's five-flow conjecture}.

\begin{con}[\cite{tut0}]\relabel{5-flow con}
Any bridgeless graph admits a nowhere-zero $5$-flow.
\end{con}

Conjecture~\ref{5-flow con} is still unproven and 
so far the best known result on this conjecture  
is due to Seymour~\cite{sey}  
who showed that every bridgelss graph has a nowhere-zero $6$-flow. For other progress regarding this 
conjecture, the reader may refer to
~\cite{jae1, jae2, maz, stef, ste}.

Seymour's result in~\cite{sey} does not imply that 
$F(G,q)>0$ holds for all bridgelss graphs $G$ and all real numbers $q\ge 6$,
although this inequality does hold for all positive integers $q\ge 6$.
This situation is similar to the case that 
the four-color Theorem 
(i.e., $P(G,4)>0$ holds for all non-loop plane graphs $G$) 
does not imply Conjecture~\ref{bir-con1}.

The study whether the inequality $F(G,q)>0$ holds
for all bridgeless graphs $G$ and real numbers $q>c$,
where $c$ is a constant,
was initiated by Welsh \cite{wel}
who proposed the following conjecture in 1970s.

\begin{con}[\cite{wel}]\relabel{welsh-con4}
For any bridgeless graph $G$, $F(G,q)>0$ holds 
for all real numbers $q\in (4,\infty)$.
\end{con}

Clearly Conjecture~\ref{bir-con1} is weaker than 
Conjecture~\ref{welsh-con4}, as it is equivalent to 
the special case of Conjecture~\ref{welsh-con4} 
when the graphs $G$ are restricted to planar graphs
and the interval for $q$ is restricted to $[4,5)$.
Counter-examples to Conjecture~\ref{welsh-con4} 
have been found while Conjecture~\ref{bir-con1} 
has neither any 
counter-example nor any result confirming 
it, even for a small interval $(a,b)$ within $(4,5)$. 

The first counter-example to Conjecture~\ref{welsh-con4} 
was due to Haggard, Pearce and Royle \cite{hag}
who showed that the generalised Petersen graph $G_{16,6}$  
has real flow roots at around $4.0252205$ and $4.2331455$, 
where the generalized Petersen graph $G_{n,k}$ 
for $n\ge 3$ and $1\le k\le \lfloor (n-1)/2\rfloor$ 
is the graph with vertex set 
$\{u_i,v_i: 1\le i\le n\}$ 
and edge set 
$\{u_iv_i, u_iu_{i+1}, v_iv_{i+k}: 1\le i\le n\}$,
$u_{n+1}$ is considered as $u_1$ and 
$v_s$ for $s>n$ is considered as $v_t$
and $t$ is the integer with $1\le t\le n$ 
such that $s-t$ is a multiple of $n$.
These graphs were introduced by Coxeter~\cite{cox} 
and named by Watkins~\cite{watk1969}. 
Clearly $G_{5,2}$ is the Petersen graph. 

There may be other counter-examples in the family of
generalized Petersen graphs to Conjecture~\ref{welsh-con4}.
But Haggard, Pearce and Royle \cite{hag} believed that
any counter-example to Conjecture~\ref{welsh-con4} 
does not have flow roots greater than or equal to $5$. 
\label{r2-q13} 
They proposed the following conjecture by   
changing the interval  to  $[5,\infty)$. 

\begin{con}[\cite{hag}]\relabel{hag-con5}
For any bridgeless graph $G$, $F(G,q)>0$ holds 
for all real numbers 
$q\in [5,\infty)$.
\end{con}

A few years ago, the above conjecture was disproved by 
Jacobsen and Salas~\cite{jaco} who found counter-examples 
by studying the 
generalised Petersen graphs $G_{6n,6}$ 
and $G_{7n,7}$ for $n\ge 2$. 
\label{r2-q14} 

\begin{theo}[\cite{jaco}]\relabel{jaco-them1}
\begin{enumerate}
\item[(i)] The value $q=5$ is an isolated accumulation point of 
real flow roots of graphs from the set 
$\{G_{6n,6}, G_{7n,7}:n\ge 3\}$. 

\item[(ii)] The value $q'\approx 5.235261$ (where $\approx $ means ``within $10^{-6}$")
is an accumulation point of real zeros of $F(G_{7n,7},q)$.
\end{enumerate}
\end{theo}

Jacobsen and Salas~\cite{jaco} further 
modified Conjecture~\ref{hag-con5}
by changing the interval to $[6,\infty)$.

\begin{con}[\cite{jaco}]\relabel{jaco-con6}
For any bridgeless graph $G$, $F(G,q)>0$ holds for all real numbers 
$q\in [6,\infty)$.
\end{con}

As pointed out by Jacobsen and Salas~\cite{jaco}, 
Conjecture~\ref{jaco-con6} might be false and 
it might even be the case that there does not exist any finite upper bound for the real flow roots of general graphs.
Now we propose the following conjecture which is much weaker 
than Conjecture~\ref{jaco-con6}.

\begin{con}\relabel{conc}
There exists a constant $c$ such that 
for any bridgeless graph $G$, $F(G,q)>0$ holds for all real numbers 
$q\ge c$.
\end{con}

Jackson~\cite{jac2} showed that for any bridgeless graph $G$
of order $n$, all real flow roots of $G$ are small than $2\log_2 n$.

\begin{theo}[\cite{jac2}]\relabel{C-th0}
For any bridgeless graph $G$ of order $n$,
$F(G,q)>0$ holds for all real $q\ge 2\log_2 n$.
\end{theo}

Theorem~\ref{C-th0} is actually  a special case of a more general
result on the characteristic polynomial 
$C(M, q)$ of a matroid $M=(E,r)$, where 
\begin{equation}\relabel{char-fun}
C(M,q)=\sum_{A\subseteq E}(-1)^{|A|}q^{r(M)-r(A)}.
\end{equation}
For any graph $G$, 
if $M_G$ and $M^*_G$ are the cycle matroid 
and the cocycle matroid of $G$ respectively, 
then $C(M_G,q)=q^{-c}P(G,q)$,  where 
$c$ is the number of components of $G$, and 
$C(M^*_G, q) = F(G,q)$.
Thus $C(G,q)$ is an extension of 
both $P(G,q)$ and $F(G,q)$.

Oxley~\cite{oxl} showed that if every cocircuit 
of $M$ has a size at most $d$,  
then $C(M, q) > 0$ holds for all real numbers $q\ge d$.
Jackson~\cite{jac2} noticed  that 
the idea in Oxley's proof can be applied to 
get a more general result. 

A {\it simple minor} of $M$ is a minor which contains no loops or circuits of length two.

\begin{theo}[\cite{jac2}]\relabel{C-th1}
Let $M$ be a matroid. 
If every simple minor of $M$ has a cocircuit of size at most $d$,
then $C(M, q) > 0$ for all real numbers $q\ge d$.
\end{theo}

As $F(G,q)=C(M^*_G,q)$, Theorem~\ref{C-th1} implies that 
for any bridgeless graph $G$,
if every 3-edge-connected minor of $G$ has a circuit of 
length at most $d$, then $F(G,q) > 0$ holds 
for all real numbers $q\ge d$. 
By Balbuena and Garc\'ia-V\'azquez's result in 
\cite{Bal2007}, 
every 3-connected graph $G$ of order $n$ 
contains a circuit of length at most $2 \log_2 n$.
Thus Theorem~\ref{C-th0} follows from Theorem~\ref{C-th1}.

\section{Graphs with real flow roots only \relabel{real zero}}

A graph $G=(V,E)$ is said to be 
a {\it chordal graph} 
if for any $S\subseteq V$ with $|S|\ge 4$, 
the subgraph of $G$
induced by $S$ is not isomorphic to any cycle. 
Dirac~\cite{dirac1961} showed that $G$ is chordal 
if and only if there is an ordering $u_1,u_2,\cdots,u_n$ 
of its vertices such that for all $i=1,2,\cdots,n$,
the subgraph of $G$ induced by 
$\{u_i\}\cup (N_G(u_i)\cap \{u_1,u_2,\cdots,u_i\})$ is a clique. By this result, 
for any chordal graph $G$, 
there exist positive integers $k_i$'s such that 
\begin{equation}\relabel{eq8-1}
P(G,q)=\prod_{i=0}^k (q-i)^{k_i},
\end{equation}
where $k=\chi(G)-1$ (see \cite{rea, read1988}).
Thus all chromatic roots of 
a chordal graph are non-negative integers. 
Some non-chordal graphs also have this property
(see \cite{dmi1980, dong3, dong4, read1975, read1988}). 

By Theorem~\ref{Tut-th1} and (\ref{eq8-1}), 
if $G$ is a plane graph and 
its dual $G^*$ is chordal, 
then all flow roots of $G$ are positive integers.
Kung and Royle \cite{kung} showed that 
the converse statement also holds. 

\begin{theo}
[\cite{kung}]
\relabel{kungroyle-theo}
If $G$ is a bridgeless graph, then 
its flow roots are integral
if and only if $G$ is the dual of a chordal and plane graph. 
\end{theo} 

Clearly, the key point in Theorem~\ref{kungroyle-theo}
is its necessity, i.e., 
\label{r2-q15} 
if $G$ has integral flow roots only,
then $G$ is the dual of a chordal and plane graph.
Does this conclusion hold when the condition 
``$G$ has integral flow roots only" is replaced by 
``$G$ has real flow roots only"?
So far it is unknown if there exists a bridgeless graph 
having real flow roots only which also contains\label{r2-q16}
non-integral real flow roots.
In \cite{dong0}, Dong studies the following problem.

\begin{prob}
\relabel{prob1} 
Is it true that any bridgeless graph 
containing real flow roots only 
has integral flow roots only?
\end{prob}

So far Problem~\ref{prob1} is even open for planar graphs.\label{r2-q17}
Let $\realr$ (resp. $\intr$) 
be the set of bridgeless graphs 
which have real (resp. integral) flow roots only.
Clearly, $\intr\subseteq \realr$.
Problem~\ref{prob1} asks if $\realr=\intr$ holds.
Dong~\cite{dong0} 
obtained the following results on the study of
problem~\ref{prob1}.

\begin{theo}[\cite{dong0}]
\relabel{main-th}
Assume that $G=(V,E)$ is a graph in $\realr$.
If some flow roots of $G$ are not in the set $\{1,2,3\}$, 
then $|V|+17\le |E|<(32|V|-49)/5$ and 
$G$ has at least $9$ flow roots in the interval $(1,2)$.
\end{theo}

Theorems~\ref{kungroyle-theo} and~\ref{main-th}
imply some equivalent statements on graphs in $\realr$.

\begin{cor}[\cite{dong0}]
\relabel{sect1-cor}
For any graph $G\in \realr$,
the following statements are equivalent:
\begin{enumerate}
\item $G$ is the dual of some plane chordal graph;
\item each flow root of $G$ is in the set 
$\{1,2,3\}$;
\item  $G$ has no flow roots in the interval $(1,2)$.
\end{enumerate}
\end{cor}

By Corollary~\ref{sect1-cor}, $\intr$ is actually
the set of bridgeless graphs $G$ whose 
flow roots are in the set $\{1,2,3\}$.

By Lemmas~\ref{block-factor} and \ref{2-edge}, to study Problem~\ref{prob1},
it suffices to consider $3$-edge connected non-separable graphs in $\realr$ which do not contain any proper $3$-edge-cut.
The article in \cite{dong0} also showed that 
if $G$ is such a graph in $\realr-\intr$,
then 
$|W(G)|\ge 3$ 
and $G$ contains at least 
$f(|W(G)|)$ flow roots in $(1,2)$, where 
$f(k)$ has values $9,11,14$ respectively for $k=3,4,5$ and  
$f(k)=\lceil \frac{27k}{11}-\frac{27}{22} \rceil$
for all integers $k\ge 6$.

So far there is no research conducted on counting 
the number of flow roots of 
a graph in the interval $(1,2)$, 
except some work of determining 
those graphs which have no real flow roots in $(1,2)$,
as mentioned in Section~\ref{zero-free}.
We end this article with the following problem,
which is related to Problem~\ref{prob1}.

\begin{prob}\relabel{prob2}
Is there a bridgeless graph $H$ 
with at least $f(k)$ flow roots in $(1,2)$,
where $k=|W(H)|\ge 3$ 
and the number of flow roots of $H$ 
counts their multiplicities?
\end{prob}

\end{document}